\documentclass[12pt]{article}
\usepackage{amssymb, amsmath}

\newtheorem{theorem}{Theorem}[section]

\newtheorem{definition}{Definition}[section]

\newcommand{\vare}{\varepsilon}
\newcommand{\n}{\nonumber}

\newcommand{\si}{\sigma_R (x)}

\newcommand{\s}{\sigma}

\newcommand{\bb}{\begin{equation}}
\newcommand{\ee}{\end{equation}}
\newcommand{\bq}{\begin{eqnarray}}
\newcommand{\eq}{\end{eqnarray}}
\newcommand{\bqn}{\begin{eqnarray*}}
\newcommand{\eqn}{\end{eqnarray*}}

\begin{document}
\title{ Liouville type of theorems for the Euler  and the Navier-Stokes equations}
\author{Dongho Chae\\
Department of Mathematics\\
              Sungkyunkwan University\\
               Suwon 440-746, Korea\\
              {\it e-mail : chae@skku.edu}}
 \date{}
\maketitle
\begin{abstract}
We prove Liouville type of theorems for weak solutions of the
Navier-Stokes and the Euler equations. In particular, if  the
pressure satisfies  $ p\in L^1 (0,T; L^1 (\Bbb R^N))$ with
$\int_{\Bbb R^N} p(x,t)dx \geq 0$, then the corresponding velocity
should be trivial, namely  $v=0$ on $\Bbb R^N \times (0,T)$. In
particular, this is the case when  $p\in L^1 (0,T; \mathcal{H}^1
(\Bbb R^N))$, where $\mathcal{H}^1 (\Bbb R^N)$ the Hardy space. On
the other hand, we have equipartition of energy over each component,
if $p\in L^1 (0,T; L^1 (\Bbb R^N))$ with $\int_{\Bbb R^N} p(x,t)dx
<0$. Similar
 results hold also for the magnetohydrodynamic
equations.
\end{abstract}

\section{Introduction}
 \setcounter{equation}{0}
We are concerned  on  the Navier-Stokes equations(the Euler
equations for $\nu=0$) on  $\Bbb R^N$, $N\in \Bbb N, N\geq 2$.
\[
\mathrm{ (NS)_\nu}
 \left\{ \aligned
 &\frac{\partial v}{\partial t} +(v\cdot \nabla )v =-\nabla p +\nu \Delta v,
 \quad (x,t)\in \Bbb R^N\times (0, \infty) \\
 & \textrm{div }\, v =0 , \quad (x,t)\in \Bbb R^N\times (0,
 \infty)\\
  &v(x,0)=v_0 (x), \quad x\in \Bbb R^N
  \endaligned
  \right.
  \]
  where $v(x,t)=(v^1 (x,t), \cdots, v^N (x,t))$ is the velocity, $p=p(x,t)$ is the
  pressure, and $\nu \geq$ is the viscosity.
  Given $a,b \in \Bbb R^N$, we denote by $a\otimes b $  the $N\times N$ matrix with $(a\otimes b)_{ij}=a_ib_j$.
  For two $N\times N$ matrices $A$ and $B$ we denote $A:B=\sum_{i,j=1}^N A_{ij} B_{ij}$.
 Given $m\in \Bbb N\cup \{0\}, q\in [1, \infty]$, we denote
 $$W^{m,q}_\s (\Bbb R^N):= \left\{ v\in [W^{m,q} (\Bbb R^N)]^N,\,\,\mathrm{ div}\, v=0 \right\},
 $$
 where $W^{m,q} (\Bbb R^N)$ is the standard Sobolev space on $\Bbb R^N,$ and
 the derivatives in  div $(\cdot)$ are in the sense of
 distribution. In particular, $H^m_\s (\Bbb R^N):=W^{m,2}_\s (\Bbb
 R^N)$ and $L^q_\s (\Bbb R^N):=W^{0, q}_\s (\Bbb R^N)$.
The Schwartz class of functions, which consists of rapidly
decreasing smooth functions, is denoted by $\mathcal{S}$ with its
dual $\mathcal{S}'$.  Let  $\varphi \in \mathcal{S} (\Bbb R^N)$ with
  $\int_{\Bbb R^N} \varphi (x) dx \neq 0$ be given. We set  $ \varphi_t (x)= t^{-N} \varphi (t^{-1} x),
  t>0$.
   Then, the  Hardy space $\mathcal{H}^q(\Bbb R^N)$, $0<q\leq 1$,  is defined
  by
$$\mathcal{H}^q(\Bbb R^N)=\left\{ f\in \mathcal{S }'\,|\,  \mathcal{M}_\varphi f(x):=\sup_{t>0} |f*\varphi_t (x) |\in L^q
  (\Bbb R^N)\right\}
  $$
 with the norm $\|f\|_{\mathcal{H}^q}:=\|\mathcal{M}_\varphi
 f\|_{L^q}$. It is well-known that the definition is independent of
 the choice of $\varphi\in \mathcal{S}$ with  $\int_{\Bbb R^N} \varphi (x) dx \neq 0$(see \cite{ste}).
 A property of $\mathcal{H}^q (\Bbb R^N)$, which will be used later
 is the fact about its dual
 \bb\label{hardy}
\left[\mathcal{H}^q (\Bbb R^N)\right]'= C^\gamma (\Bbb R^N), \quad
\gamma=N\left(\frac{1}{p} -1\right),
 \ee
 if $0<p<1$, where $C^\gamma (\Bbb R^N)$ is the homogeneous
 H\"{o}lder space.
  In $\Bbb R^N$ we define  weak
 solutions of the Navier-Stokes(Euler)
equations as follows.
  \begin{definition} We say the pair $(v,p)\in L^1 (0, T; L^2_\s (\Bbb R^N))\times L^1
  (0, T; \mathcal{S}' (\Bbb R^N ))$ is a weak solution of $(NS)_\nu$
  on $ \Bbb R^N \times(0, T )$
  if
  \bq\label{11}
  \lefteqn{-\int_{0}^{T} \int_{\Bbb R^N} v(x,t) \cdot \phi(x)\xi ' (t)dxdt -\int_{0}^{T}\int_{\Bbb R^N} v(x,t)
  \otimes v(x,t):\nabla \phi (x) \xi (t) dxdt}\hspace{.in} \n \\
  && = \int_{0} ^{T}< p(t),\mathrm{ div }\, \phi >
  \xi (t)dt +\nu \int_{0}^{T} \int_{\Bbb R^N} v(x,t)\cdot \Delta \phi (x) \xi (t)dxdt\n \\
  \eq
   for all $\xi \in C^1 _0 (0, T)$ and $\phi =[C_0 ^\infty (\Bbb R^N )]^N
   $, where $<\cdot,\cdot >$ denotes the dual pairing between
   $\mathcal{S}$ and $\mathcal{S}'$.
  \end{definition}
  The definition is weaker than the standard Leray-Hopf weak solution for the
Navier-Stokes equations, since
 we are concerned also on possible weak solutions of the Euler
 equations, the right function space of whose existence is not yet known. Below we denote
$$\mathcal{E}_j (t)=\frac12 \int_{\Bbb R^N} (v^j (x,t))^2 dx, \quad
j=1, \cdots, N
$$
which will be called the $j-$th component of the total energy,
$$\mathcal{E}(t)=\frac12 \int_{\Bbb R^N} |v (x,t)|^2 dx=\mathcal{E}_1 (t)+\cdots+\mathcal{E}_N (t).
$$
Let us introduce the function class,
  $$L^1_\pm (0,T; L^1(\Bbb R^N))=\left\{ f\in L^1 (0,T; L^1(\Bbb R^N)),\, \int_{\Bbb R^N} f(x,t)dx \geq(\leq)\, 0 \,\mathrm{ a.e.} \, t\in (0, T)\right\}.
  $$
 \begin{theorem} Let  $(v,p)$ be a weak solution to $(NS)_\nu$
 with $\nu \geq 0$.
\begin{itemize}
\item[(i)] (Liouville type of property) Suppose
  \bb\label{19}
  \mbox{either}\qquad p\in L^1 _+ (0, T;L^1 (\Bbb R^N
  )),\quad\mbox{or}\quad p\in L^1 (0, T; \mathcal{H}^q (\Bbb R^N))
  \ee
 for some $q\in (0, 1]$. Then,
 \bb
 \label{112a} \mbox{$v(x,t)=0$  a.e.  in $ \Bbb R^N \times
(0,T)$},
 \ee
\item[(i)] (Equipartition of energy) Suppose $p\in L^1 _-(0, T;L^1 (\Bbb R^N
  ))$. Then,
 \bb\label{111} \mathcal{E}_1 (t)=\cdots=\mathcal{E}_N
(t)=-\frac12\int_{\Bbb R^N} p(x,t)dx,
 \ee
  and
 \bb\label{112}
\int_{\Bbb R^N} v^j (x,t)v^k(x,t)dx=0\quad \forall j,k\in
\{1,\cdots, N\} \quad\mbox{with}\quad j\neq k
 \ee
  for almost every
$t\in (0, T)$.
\end{itemize}
\end{theorem}
{\it Remark 1.1 } Let us recall that  $\int_{\Bbb R^N} f(x)dx =0$,
if $f\in \mathcal{H}^1   (\Bbb R^N )$, where $\mathcal{H}^1   (\Bbb
R^N )$ is the Hardy space in $\Bbb R^N$(see \cite{ste}), and
 $$ L^1(0,T; \mathcal{H}^1   (\Bbb R^N ))\subset L^1_+ (0,T; L^1(\Bbb R^N)).
 $$
The part (i) of the above theorem says that the cancelation property
of the pressure is nt allowed for nontrivial solutions of the
Navier-Stokes and the Euler equations.  Note that the condition
(\ref{19}) with $q=1$ is already far beyond the natural scaling of
the usual regularity criterion on the pressure for the Navier-Stokes
equations,
$$
p\in L^q (0, T; L^r (\Bbb R^N )), \qquad \frac{2}{q}
+\frac{N}{r}\leq 2.
$$
(see \cite{cha, ber}), and our conclusion is not just the solution
is
regular, but it is {\em trivially zero.}\\
\ \\
\noindent{\it Remark 1.2}
 We also recall the relation
between the pressure and velocity for the Navier-Stokes and Euler
equations:
 \bb\label{18}
p(x,t)=\sum_{j,k=1}^N \left[R_jR_k( v^j (\cdot,t) v^k(\cdot,t
)\right] (x),
 \ee
  where $R_j, j=1, \cdots, N,$ is the Riesz transforms in $\Bbb
  R^N$, defined by
 $$R_j (f)(x)= C_N \lim_{\vare\to 0} \int_{|y|>\vare} \frac{y_j}{|y|^{n+1}} f(x-y)
 dy, \quad C_N =\frac{\Gamma \left(\frac{N+1}{2}\right)}{\pi
 ^{(N+1)/2}}
 $$
Thus, we  find that (\ref{19}) with $q=1$  is guaranteed if
  \bb\label{hardy}
v\otimes v \in L^1(0,T;\mathcal{H}^1 (\Bbb R^N) ).
 \ee
   In reality the
pressure for the Leray weak solutions of the N-dimensional ($N=2,3$)
Navier-Stokes equations has the property that $v\in L^2(0,T;
H^1_\s(\Bbb R^N))\cap L^\infty (0, T; L^2_\s (\Bbb R^N
))$(\cite{ler}), which implies
 \bb\label{weak}
v\otimes v\in L^1(0,T; L^q (\Bbb R^N))
 \ee
   for all $q\in [1,
\frac{N}{N-2}]$ if $N\geq 3$, while $q\in [1, \infty)$ if $N=2$. On
the other hand, we note that the local smooth solution $v\in C([0,
T); H^m _\s (\Bbb R^N))$, $m>N/2+1$, constructed by
Kato(\cite{kat}), has the property $v \in C([0, T); L^q_\s (\Bbb
R^N))$ for all $q\in [2, \infty]$ due to the embedding $H^m (\Bbb
R^N)\hookrightarrow L^\infty (\Bbb R^N)$, combined with the
interpolation between $L^2(\Bbb R^N)$ and $L^\infty(\Bbb R^N)$.
Hence, for $v\in C([0, T); H^m _\s (\Bbb R^N))$, $m>N/2+1$, we have
\bb \label{ka}
 v\otimes v \in C([0,T);L^q(\Bbb R^N) )\quad \forall
q\in [1, \infty].
 \ee
 It would be interesting to recall  the related
known properties of the pressure for the Leray weak solutions, which
are proved in \cite{coi} (see also \cite{lio}):
 \bqn
D^2 p &\in& L^1 (0, T; \mathcal{H}^1   (\Bbb R^N )),\\
\nabla p &\in & L^2 (0, T; \mathcal{H}^1   (\Bbb R^N ))\cap L^1(0,T; L^{\frac{N}{N-1}, 1} (\Bbb R^N)),\\
p &\in& \left\{ \aligned & L^1(0,T; L^{\frac{N}{N-2}, 1} (\Bbb R^N)), \, N\geq 3\\
                      &L^1(0,T; \mathcal{C}_0 (\Bbb R^2)), \, N=2\endaligned \right.
\eqn where $L^{q,r} (\Bbb R^N)$ is the Lorentz space, and
$\mathcal{C}_0 (\Bbb R^2)$ is the class of continuous functions
vanishing near
infinity\\
\ \\
\noindent{\it Remark 1.3} One immediate corollary of the above
theorem is that the following pressureless Navier-Stokes(Euler)
system is locally ill-posed if $v_0 \in H^m_\s (\Bbb R^N)$,
$m>N/2+1$,
\[
 \left\{ \aligned
 &\frac{\partial v}{\partial t} +(v\cdot \nabla )v = \nu \Delta v,
 \quad (x,t)\in \Bbb R^N\times (0, \infty) \\
 & \textrm{div }\, v =0 , \quad (x,t)\in \Bbb R^N\times (0,
 \infty)\\
  &v(x,0)=v_0 (x), \, \textrm{div }\, v_0 =0,\quad x\in \Bbb R^N
  \endaligned
  \right.
  \]
   since we need to have $v(\cdot, t)=0$
for $t>0$ from the fact $\int_{\Bbb R^N} p(x,t)dx =0$ for all
$t>0$.\\
\ \\
 \noindent{\it Remark 1.4}
 In \cite{jiu} the equipartition of energy over each component
 has been derived for {\em steady} Euler equations in a different
 context, using different definition of weak solutions. A completely different type of approach to the
Liouville type of theorems for the Navier-Stokes equations is
studied in \cite{koc}.\\

\section{Proof of the Main Theorem}
 \setcounter{equation}{0}

\noindent{\bf Proof of Theorem 1.1}\\
\noindent{\underline{(a) the case $p\in L^1 (0, T; L^1 (\Bbb R^N))$
}:}
 Let us
consider a cut-off function $\sigma\in C_0 ^\infty(\Bbb R^N)$ such
that
$$
\sigma (x)= \sigma(|x|)=\left\{ \aligned
                  &1 \quad\mbox{if $|x|<1$}\\
                     &0 \quad\mbox{if $|x|>2$},
                      \endaligned \right.
$$
and $0\leq \sigma  (x)\leq 1$ for $1<|x|<2$.
Then, given $R >0$, we set
$$\varphi_R (x)=\frac{x_1^2}{2}\sigma
\left(\frac{x}{R} \right). $$
 Let $\xi \in C^1 _0 (0, T)$. We  choose the vector test function
 $\phi$ in (\ref{11}) as
$$\phi= \nabla \varphi_R (x)=\left(x_1\sigma_R (x) + \frac{x_1 ^2}{2} \partial_1 \si , \frac{x_1^2}{2}
\partial_2 \si , \cdots , \frac{x_1^2}{2} \partial_N \si \right).
$$
 Then, (\ref{11}) becomes
 \bq\label{21}
\lefteqn{0=\int_0 ^T\int_{\Bbb R^N} (v^1 (x,t))^2\sigma_{R}(x) \xi(t)  \,dx dt}\n\\
&&\quad+\int_0 ^T\int_{\Bbb R^N}  p (x,t) \sigma_R (x)\xi(t)\,dx dt\n \\
&&\quad+\int_0 ^T\int_{\Bbb R^N} (v^1 (x,t))^2\left[2 x_1
\partial_1\sigma_{R}(x) +\frac{x_1^2}{2} \partial_1 ^2 \si \right]\xi(t) \,dx
dt\n \\
 &&\quad+2\sum_{j=2}^N \int_0 ^T\int_{\Bbb R^N} v^1 (x,t)v^j
(x,t)\left[ x_1\partial_j  \sigma_{R}(x) +\frac{x_1^2}{2} \partial_1
\partial_j
\si \right]\xi(t)\,dxdt\n \\
&&\quad+\sum_{j,k=2}^N\int_0 ^T\int_{\Bbb R^N}  v^j(x,t)v^k (x,t)
x_1^2\partial_j \partial_k \sigma_R (x)
 \xi(t) \,dxdt\n \\
 &&\quad+ \int_0 ^T\int_{\Bbb R^N}p(x,t) \left[ 2x_1 \partial_1 \si +\frac{x_1^2}{2}
 \Delta \sigma_R (x)\right] \xi(t) \, dxdt\n \\
 &&:=I_1+\cdots +I_6.
\eq
 Note that the first term of the left hand side  and the second term of the right hand side in (\ref{11})
 vanish, since
 $$
 \int_0 ^T\int_{\Bbb R^N} v(x,t)\cdot \nabla\varphi_R (x) \xi'(t) dxdt=0,
 $$
 and
$$
 \int_0 ^T\int_{\Bbb R^N} v(x,t)\cdot \nabla (\Delta \varphi_R (x)) \xi(t) dxdt=0
 $$
 for  $v\in L^1 (0, T; L^2_\s (\Bbb R^N))$ by the definition
 of divergence free condition in the sense of distribution.
 We pass $R\to \infty $ in (\ref{21}).
Since $v\in L^\infty(0, T; L^2 _\s (\Bbb R^N))$ by hypothesis,
 \bq\label{dom}
\lefteqn{\left|I_1-\int_0 ^T\int_{\Bbb R^N} (v_1 (x,t))^2
\xi(t)\,dxdt\right|
  \leq \int_0 ^T\int_{\Bbb R^N} (v_1 (x,t))^2 |\xi(t)| |1-\si
 |\,dxdt}\hspace{1.5in}\n \\
&& \leq \sup_{0<t<T} |\xi(t)| \int_0 ^T \int_{|x|>R}
 (v^1(x,t))^2 dxdt
 \eq
 Since
 $$g_R(t):=\int_{|x|>R} (v_1 (x,t))^2 \,dx \to 0 \quad \mbox{as}\quad R\to \infty
$$
for almost every $t\in (0, T)$, and
$$ |g_R (t)|\leq g(t):=\int_{\Bbb R^N} (v_1 (x,t))^2 \,dx
$$
with $g\in L^1 (0,T)$, we can apply the dominated convergence
theorem in (\ref{dom}) to get
$$
\int_0 ^T \int_{|x|>R}
 (v^1(x,t))^2 dxdt \to 0\quad \mbox{as}\quad R\to \infty,
 $$
 and hence
 \bb\label{22}
 I_1\to\int_0 ^T\int_{\Bbb R^N} (v_1 (x,t))^2 \xi(t)\,dxdt \quad
 \mbox{as$\quad R\to \infty$}.
 \ee
Since $p\in L^1 (0, T; L^1 (\Bbb R^N))$, we have
\bqn
\left| I_2-\int_0 ^T\int_{\Bbb R^N}  p (x,t) \xi(t)\,dx dt\right|&\leq&\sup_{0<t<T} |\xi(t)|\int_0 ^T\int_{\Bbb R^N}  |p (x,t)||1-\si |\,dx dt\n \\
& \leq& \sup_{0<t<T} |\xi(t)| \int_0 ^T \int_{|x|>R} |p(x,t)| dx dt
\to 0 \eqn as $R\to \infty$ by the dominated convergence theorem.
Hence,
 \bb\label{22a}
 I_2 \to\int_0 ^T\int_{\Bbb R^N}  p (x,t) \xi(t)\,dx dt \quad \mbox{as $R\to \infty$}
 \ee
 We will show below that $I_3, \cdots, I_6 \to 0$ as $R\to \infty$.
In view of (\ref{hardy})
 Next, we note that for $m\geq 1$ and $j,k \in \{ 1, \cdots, N\}$
\bq\label{23}
 \lefteqn{\left|\int_0 ^T \int_{\Bbb R^N} \xi (t) v^j(x,t)v^k (x,t) x_1 ^m
D^m \si dxdt\right|}\hspace{.0in} \n \\
&&\leq\frac{1}{R^m} \sup_{1<s<2}|\sigma ^{(m)} (s)| \int_0 ^T
\int_{R<|x|<2R} |\xi (t)| |v (x,t)|^2
|x|^m dxdt\n \\
 &&\leq 2^{m}\sup_{1<s<2}|\sigma ^{(m)} (s)| \sup_{0<t<T} |\xi(t)|\int_0
^T\int_{R<|x|<2R} |v(x,t)|^2 dxdt \n \\
&&\to 0
 \eq
 as $R\to \infty$ by the dominated convergence theorem, which shows
 that $I_3, I_4$ and $I_5$ converge to zero as $R\to \infty$.
Similarly,
  \bq\label{24}
 \lefteqn{\left|\int_0 ^T \int_{\Bbb R^N} \xi (t)p(x,t) x_1 ^m
D^m \si dxdt\right|}\hspace{.0in} \n \\
&&\leq\frac{1}{R^m} \sup_{1<s<2}|\sigma ^{(m)} (s)| \int_0 ^T
\int_{R<|x|<2R} |\xi (t)| |p (x,t)|
|x|^m dxdt\n \\
 &&\leq 2^{m}\sup_{1<s<2}|\sigma ^{(m)} (s)| \sup_{0<t<T} |\xi(t)|\int_0
^T\int_{R<|x|<2R} |p(x,t)| dxdt \n \\
&&\to 0
 \eq
 as $R\to \infty$ by the dominated convergence theorem, which shows
 that $I_6$ converges to zero as $R\to \infty$, since $p\in
  L^1 (0, T; L^1 (\Bbb R^N ))$.
Therefore, after passing $R\to \infty$ in (\ref{21}), we are left
with
 $$
 \int_0 ^T \int_{\Bbb R^N} \xi(t) \left[(v^1 (x,t))^2+p(x,t)\right]dxdt
 =0 \quad \forall \xi \in C^1_0 (0,T).
 $$
 Hence,
 $$\mathcal{E}_1 (t)=-\frac12\int_{\Bbb R^N} p(x,t)dx \quad \mbox{for
 almost every $t\in (0,T)$}.
 $$
 Similarly, if we choose the vector test function $\phi$ in (\ref{11}) as
 $$\phi(x)=\nabla
 \left(\frac{x_j^2}{2} \si \right), \quad j\in \{1, \cdots, N\}
 $$
 then we could obtain
 $$
 \int_0 ^T \int_{\Bbb R^N} \xi(t) \left[(v^j (x,t))^2+p(x,t)\right]dxdt
 =0 \quad \forall \xi \in C^1_0 (0,T),
 $$
 and hence
$$\mathcal{E}_j (t)=-\frac12\int_{\Bbb R^N} p(x,t)dx \quad \mbox{for
 almost every $t\in (0,T)$}
 $$
 for all $j\in \{ 1, \cdots, N\}$. This proves (\ref{111}).
 In order to prove (\ref{112}) we choose the test function
\bqn
\lefteqn{\phi(x)=\nabla \left(x_1 x_2 \si \right)}\\
&&=(x_2 \si +x_1x_2 \partial_1 \si, x_1 \si +x_1x_2 \partial_2 \si ,x_1x_2
  \partial_3\si ,\n \\
  &&\hspace{1.in} \cdots, x_1x_2 \partial_N \si ).
\eqn
 Then, we have
\bq\label{26}
 \lefteqn{ 0=2\int_0 ^T \int_{\Bbb R^N} v^1 (x,t)v^2(x,t) \si\xi(t)\,
 dxdt}\n \\
 &&\quad+\int_0 ^T \int_{\Bbb R^N} (v^1(x,t) )^2 \left[ 2x_2
 \partial_1 \si + x_1x_2 \partial_1^2 \si \right] \xi(t) \,dxdt \n \\
 &&\quad +\int_0 ^T \int_{\Bbb R^N} (v^2(x,t) )^2 \left[ 2x_1
 \partial_1 \si + x_1x_2 \partial_2^2 \si \right] \xi(t)\, dxdt\n \\
 &&\quad +\int_0 ^T \int_{\Bbb R^N} v^1 (x,t)v^2(x,t)\left[ x_1
 \partial_1 \si +x_2
\partial_2 \si  +2 x_1x_2 \partial_1\partial_2 \si \right] \xi(t)\,
 dxdt\n \\
  && \quad+ 2\sum_{j=3}^N\int_0 ^T \int_{\Bbb R^N} v^1
(x,t)v^j(x,t)
 \left[ x_2 \partial_j \si +x_1x_2 \partial_1 \partial_j \si \right]
 \xi(t)\, dxdt \n \\
 && \quad+ 2\sum_{j=3}^N\int_0 ^T \int_{\Bbb R^N} v^2 (x,t)v^j(x,t)
 \left[ x_1 \partial_j \si +x_1x_2 \partial_1 \partial_j \si \right]
 \xi(t)\, dxdt\n \\
&& \quad+ 2\sum_{j,k=3}^N\int_0 ^T \int_{\Bbb R^N} v^j
(x,t)v^k(x,t)x_1x_2 \partial_j \partial_k \si
 \xi(t)\, dxdt\n \\
 &&\quad+ \int_0 ^T \int_{\Bbb R^N} p(x,t)\left[ 2x_2 \partial_1 \si
 + 2x_1\partial_2\si + x_1x_2 \Delta \si \right] \xi(t) dxdt \n \\
  &&:=J_1+\cdots +J_{8}.
  \eq
  Similarly to the previous proof we deduce
 \bqn
 \lefteqn{\left|J_1-2\int_0 ^T \int_{\Bbb R^N} v^1 (x,t)v^2(x,t) \xi(t)\,
 dxdt\right|}\n \\
 &&\leq 2\int_0 ^T \int_{\Bbb R^N}|v^1 (x,t)v^2(x,t)|\,|\xi(t)|\, |1-\si
 |\,dxdt\n \\
 &&\leq 2\sup_{0<t<T} |\xi(t)|\int_0 ^T \int_{|x|>R } |v(x,t)|^2 dx
 dt\to 0
 \eqn
 by the dominated convergence theorem, and
 $$ J_1 \to 2\int_0 ^T \int_{\Bbb R^N} v^1 (x,t)v^2(x,t) \xi(t)\,
 dxdt.$$
  Using the computations similar to (\ref{23}) and (\ref{24}), we also find that
  $$ \sum_{k=2}^{13}|J_k| \to 0 \qquad\mbox{as $R\to \infty$}. $$
  Hence, taking $R\to \infty$ in (\ref{26}), we are left with
  $$0=\int_0 ^T \int_{\Bbb R^N} v^1 (x,t)v^2(x,t) \xi(t)\,
 dxdt\quad\forall \xi \in C^1_0 (0,T),
 $$
 and therefore
 $$\int_{\Bbb R^N} v^1(x,t)v^2(x,t)dx =0 \quad\mbox{for almost every
 $t\in (0,T)$}.
 $$
 Similarly, choosing the test function
$$\phi(x)=\nabla (x_j x_k \si ), \quad j,k \in \{1,\cdots, N\}, j\neq k $$
and repeating the above argument, we could drive
$$\int_{\Bbb R^N} v^j(x,t)v^k(x,t)dx =0 \quad\mbox{for almost every
 $t\in (0,T)$}.
 $$
 for all $j,k \in \{1,\cdots, N\}$ with $ j\neq k$.\\
 \ \\
\noindent{\underline{(b) the case $p\in L^1 (0, T; \mathcal{H}^q
(\Bbb R^N)), 0<q\leq 1$ }:}  The borderline case $p=1$ is contained
in the part (a) above(see Remark 1.1), and we assume here $0<p<1$.
 In order to derive the Liouville type of property in this case it
 suffice to show that $I_2, I_6 \to 0$ as $R\to \infty$ in
 (\ref{21}). This can be shown by the following estimates for
 $m\in \Bbb N \cup\{0\}$
 \bqn
  \lefteqn{\left|\int_0^T \xi (t) <p (\cdot, t), x_1 ^m D^m \si
 >dt\right|  }\hspace{.5in}\n \\
 &&\leq \sup_{0\leq t\leq T} \int_0 ^T |\xi
(t)|\|p(t)\|_{\mathcal{H}^q} dt  \|x_1 ^m D^m \si
\|_{C^{\gamma}}\n\\
&&\leq \frac{C}{ R^\gamma}  \|\xi \|_{L^\infty (0, T)} \|p(t)\|_{L^1
(0, T; \mathcal{H}^q)}\to 0,
 \eqn
  where $ \gamma=  N(1/q -1 ) >0$, where we used the duality, $[\mathcal{H}^q
  (\Bbb R^N)]'= C^\gamma (\Bbb R^N)$, and the  simple estimate,
  $$ \|x^m_1 D^m \si \|_{C^\gamma}\leq  \frac{C}{ R^\gamma} .$$
$\square$\\
\ \\
\noindent{\it Remark  after the proof } A natural question is if
there exists an initial data $v_0\in H^m_\s(\Bbb R^N)$, $m>N/2+1$,
with div $v_0 =0$ such that the corresponding initial pressure
satisfies
  \bb\label{initial}
p_0=\sum_{j,k=1}^N R_jR_k (v^j_0 v^k_0) \in L^1 (\Bbb R^N),
 \ee
 but $\mathcal{E}_j(0) \neq\mathcal{E}_k(0)$ for some $j,k\in \{1,
\cdots, N\}$, $j\neq k$. If this is possible, then it implies that
function $t\mapsto \|p(t)\|_{L^1}$ is discontinuous at $t=0+$ for
the local classical solution $(v(\cdot ,t), p(\cdot ,t))$
constructed by Kato(\cite{kat}) with such initial data. Using the
Fourier transform, Professor P. Constantin showed that there exists
no such initial data(\cite{con}). Actually similar conclusion can be
derived by slight change of the above proof as follows. From the
relation (\ref{initial}) we have
$$-\int_{\Bbb R^N} p_0 (x)\Delta \psi dx=\sum_{\j,k=1}^N\int_{\Bbb
R^N} v_0 ^j (x) v_0 ^k (x) \partial_j \partial_k \psi (x)dx\,\,
\forall \psi \in C_0^2 (\Bbb R^N).
$$
Similarly to the above proof, choosing $\psi (x)=x_j x_k \si $, and
then passing $R\to \infty$, we have
 $$-\int_{\Bbb R^N} p_0 (x)dx=\int_{\Bbb R^N} (v_0 ^j (x))^2 dx
 \quad \forall j=1, \cdots, N$$
 for $j=k$, while
 $$\int_{\Bbb R^N} v_0 ^j (x)v_0 ^k (x)dx=0 \quad\forall j,k=1,
 \cdots, N,
 $$
 for $j\neq k$.\\

\section{Remarks on the MHD equations}
 \setcounter{equation}{0}

In this section we extend the previous results on the system
$(NS)_\nu$ to the magnetohydrodynamic equations in $\Bbb R^N$,
$N\geq 2$.
\[
\mathrm{ (MHD)_{\mu,\nu} }
 \left\{ \aligned
 &\frac{\partial v}{\partial t} +(v\cdot \nabla )v =(b\cdot\nabla)b-\nabla (p +\frac12 |b|^2) +\nu \Delta v, \\
 &\frac{\partial b}{\partial t} +(v\cdot \nabla )b =(b \cdot \nabla
 )v+\mu \Delta b,\\
 &\quad \textrm{div }\, v =\textrm{div }\, b= 0 ,\\
  &v(x,0)=v_0 (x), \quad b(x,0)=b_0 (x)
  \endaligned
  \right.
  \]
where $v=(v_1, \cdots , v_N )$, $v_j =v_j (x, t)$, $j=1,\cdots,N$,
is the velocity of the flow, $p=p(x,t)$ is the scalar pressure,
$b=(b_1, \cdots , b_N )$, $b_j =b_j (x, t)$, is the magnetic field,
and $v_0$, $b_0$ are the given initial velocity and magnetic field,
 satisfying div $v_0 =\mathrm{div}\, b_0= 0$, respectively. Let us
 begin with the definition of the weak solutions of $(MHD)_{\mu,\nu}$.
   \begin{definition} We say the triple $(v,b,p)\in [L^1 (0, T; L^2_\s (\Bbb R^N))]^2\times L^1 (0, T; \mathcal{S}' (\Bbb R^N ))$ is a weak solution of $ (MHD)_{\mu,\nu}$
  on $ \Bbb R^N \times(0, T )$,
  if
  \bq\label{31}
  \lefteqn{-\int_{0}^{T} \int_{\Bbb R^N} v(x,t) \cdot \phi(x)\xi ' (t)dxdt -\int_{0}^{T}\int_{\Bbb R^N} v(x,t)
  \otimes v(x,t):\nabla \phi (x) \xi (t) dxdt} \n \\
  &&=-\int_{0}^{T}\int_{\Bbb R^N}b(x,t)
  \otimes b(x,t):\nabla \phi (x) \xi (t) dxdt
 +\int_{0} ^{T} <p (t),\mathrm{ div }\, \phi >
  \xi (t)dt \n \\
  &&\quad+\frac12 \int_{0} ^{T} \int_{\Bbb R^N}|b(x,t)|^2\,\mathrm{ div }\, \phi (x)
  \xi (t)dxdt
  +\nu \int_{0}^{T} \int_{\Bbb R^N} v(x,t)\cdot \Delta \phi (x) \xi (t)dxdt,\n \\
  \eq
  and
  \bq\label{32}
  \lefteqn{-\int_{0}^{T} \int_{\Bbb R^N} b(x,t) \cdot \phi(x)\xi ' (t)dxdt -\int_{0}^{T}\int_{\Bbb R^N} v(x,t)
  \otimes b(x,t):\nabla \phi (x) \xi (t) dxdt} \n \\
  &&=-\int_{0}^{T}\int_{\Bbb R^N}b(x,t)
  \otimes v(x,t):\nabla \phi (x) \xi (t) dxdt+\mu \int_{0}^{T} \int_{\Bbb R^N} b(x,t)\cdot \Delta \phi (x) \xi (t)dxdt\n \\
  \eq
   for all $\xi \in C^1 _0 (0, T)$ and $\phi =[C_0 ^\infty (\Bbb R^N )]^N
   $.
  \end{definition}

 We have the following theorem.
\begin{theorem}
Suppose  $(v,b,p)$ ia a weak solution to $(MHD)_{\mu, \nu}$  with
$\mu, \nu \geq 0$ on $\Bbb R^N\times (0, T )$ satisfying
  \bb\label{33}
   \mbox{either}\qquad p\in L^1 _+ (0, T;L^1 (\Bbb R^N
  )),\quad\mbox{or}\quad p\in L^1 (0, T; \mathcal{H}^q (\Bbb R^N))
  \ee
 for some $q\in (0, 1]$.
Then,  for $N\geq 3$, we have
 \bb\label{35}
  v(x,t)=b(x,t)=0 \quad\mbox{a.e. in $\Bbb R^N\times(0,T)$},
 \ee
 while for $N=2$
 \bb\label{36}
 v(x,t)=0, \quad\mbox{a.e. in $\Bbb R^2\times(0,T)$},
 \ee
 \bb
  \label{37}
 \int_{\Bbb R^2} (b^1 (x,t))^2 dx= \int_{\Bbb R^2} (b^2
 (x,t))^2 dx
 \ee
 almost everywhere in $(0, T)$, and $b (x,t)$ is a
 weak solution of the heat equation $\partial_t b =\mu \Delta b
 $.
\end{theorem}
\noindent{\bf Proof  } The method of proof is similar to that of
Theorem 1.1 with slight changes. We will be brief, describing only
essential points. We choose the vector test function $\phi =\nabla
(\frac{x_j^2}{2} \si )$ with $j\in \{ 1, \cdots, N\}$ in (\ref{31}).
Then, in the case $p\in L^1 (0, T; L^1 (\Bbb R^N))$, we obtain that
 \bq\label{38}
 \lefteqn{-\int_0 ^T\int_{\Bbb R^N} (v^j(x,t))^2 \si \xi(t)dxdt=-\int_0 ^T\int_{\Bbb R^N} (b^j(x,t))^2 \si
 \xi(t)dxdt}\hspace{.0in}\n \\
 && +\frac12\int_0 ^T \int_{\Bbb R^N} |b(x,t)|^2\si \xi
 (t)dxdt + \int_0 ^T  \int_{\Bbb R^N} p(x,t)dx \si\xi (t)dxdt+  o(1),\n \\
 \eq
where $o(1)$ denotes the sum of the terms vanishing as $R\to
\infty$. Taking $R\to \infty$ in (\ref{38}), and summing over $j=1,
\cdots, N$, we find that
 \bq\label{39}
 \lefteqn{-\int_0 ^T\int_{\Bbb R^N} |v(x,t)|^2  \xi(t)dxdt-\frac{N-2}{2}\int_0 ^T\int_{\Bbb R^N} |b(x,t)|^2
  \xi(t)dxdt}\hspace{1.5in}\n \\
  &&=N\int_0 ^T  \int_{\Bbb R^N} p(x,t)dx \xi (t)dxdt .\n \\
  \eq
If $p\in L^1_+ (0, T; L^1 (\Bbb R^N))$,  choosing $\xi \in C^1_0
(0,T)$ with $\xi(t)\geq 0$ in
  (\ref{39}), then $N\geq 3$ implies $v(x,t)=b(x,t)=0$ and  $p(x,t)=0$ for almost
  every $(x,t)\in \Bbb R^N\times (0, T)$. If $N=2$, instead, then
  $v(x,t)=0$ and $p(x,t)=0$ for almost
  every $(x,t)\in \Bbb R^N\times (0, T)$. Then, (\ref{32}) becomes
  $$
-\int_{0}^{T} \int_{\Bbb R^N} b(x,t) \cdot \phi(x)\xi ' (t)dxdt=\mu \int_{0}^{T} \int_{\Bbb R^N} b(x,t)\cdot \Delta \phi (x) \xi (t)dxdt
\quad
$$
for all $\phi \in [C_0 ^2 (\Bbb R^N)]^N$ and $\xi \in C^1 _0 (0, T)$, which shows that $b(x,t)$ is  a weak solution of $\partial_t b=\mu \Delta b$.  Moreover, after passing $R\to \infty$, the equation (\ref{38}) with $j=1$ becomes
  $$\int_0 ^T\int_{\Bbb R^2}\left[ (b^2 (x,t))^2
  -(b^1 (x,t))^2 \right] \xi(t) dxdt=0,\quad\forall \xi\in C_0 ^1 (0,T)
  $$
and hence we obtain (\ref{37}). In the case $p\in L^1 (0, T;
\mathcal{H}^q (\Bbb R^N))$, $0<q<1$, following the same argument as
in the  part (b) of proof of Theorem 1.1  in the previous section,
one can derive
 $$
-\int_0 ^T\int_{\Bbb R^N} |v(x,t)|^2  \xi(t)dxdt-\frac{N-2}{2}\int_0
^T\int_{\Bbb R^N} |b(x,t)|^2
  \xi(t)dxdt=0
  $$
  instead of (\ref{39}), from which our conclusion follows.
  $\square$\\
  \[\mbox{\bf Acknowledgements}\]
  The author would like to thank to  P. Constantin and J. Verdera for useful discussion,
  and to J. Lee for reading the
  manuscript. This work was supported partially by  KRF Grant(MOEHRD, Basic
Research Promotion Fund).  Part of the research was done, while the
author was visiting University of Chicago.

\end{document}